\documentclass[a4paper,11pt]{article} 
 
\usepackage{mathrsfs}
\usepackage{amsmath} 
\usepackage[applemac]{inputenc}
\usepackage{amsfonts}
\usepackage{amssymb}
\usepackage{amsthm}
\usepackage{stmaryrd}
\usepackage{color}

\newcommand{\N}{\mathbb{N}}

\newcommand{\R}{\mathbb{R}}

\def\1{\mbox{1\hspace{-0.25em}l}}

\def\bx{{\mathbf{x}}}
\def\x{{\boldsymbol{x}}}
\def\p{{\mathbf{p}}}

\def\l{\lambda}
\def\<{\langle}
\def\>{\rangle}

\def\x{{\xi}}
\def\z{{\zeta}}

\numberwithin{equation}{section}

\newcommand{\be}{\begin{eqnarray}}
\newcommand{\ee}{\end{eqnarray}}
\newcommand{\ce}{\begin{eqnarray*}}
\newcommand{\de}{\end{eqnarray*}}
\newtheorem{theorem}{Theorem}[section]
\newtheorem{lemma}[theorem]{Lemma}
\newtheorem{remark}[theorem]{Remark}
\newtheorem{definition}[theorem]{Definition}
\newtheorem{proposition}[theorem]{Proposition}
\newtheorem{Examples}[theorem]{Example}
\newtheorem{corollary}[theorem]{Corollary}
\newtheorem{Assumption}[theorem]{Assumption}

\def\eps{\varepsilon}

\def\a{\alpha}

\def\p{\partial}
\def\d{\delta}

\def\l{\lambda}

\def\[{{\Big[}}
\def\]{{\Big]}}
\def\<{{\langle}}
\def\>{{\rangle}}
\def\({{\big(}}
\def\){{\big)}}

\def\bx{{\mathbf{x}}}

\def\bb2{{\boldsymbol{2}}}

\def\={&\!\!=\!\!&}

\def\cD{{\mathcal D}}

\def\cG{{\mathcal G}}

\def\1{{\mathbf{1}}}

\def\geq{\geqslant}
\def\leq{\leqslant}
\def\ge{\geqslant}
\def\le{\leqslant}

\def\eps{\varepsilon}

\def\a{\alpha}

\def\p{\partial}

\def\l{\lambda}

\def\[{{\Big[}}
\def\]{{\Big]}}
\def\<{{\langle}}
\def\>{{\rangle}}

\def\bx{{\mathbf{x}}}

\def\={&\!\!=\!\!&}
\def\bt{\begin{theorem}}
\def\et{\end{theorem}}
\def\bl{\begin{lemma}}
\def\el{\end{lemma}}
\def\br{\begin{remark}}
\def\er{\end{remark}}
\def\bx{\begin{Examples}}
\def\ex{\end{Examples}}
\def\bd{\begin{definition}}
\def\ed{\end{definition}}
\def\bp{\begin{proposition}}
\def\ep{\end{proposition}}
\def\bc{\begin{corollary}}
\def\ec{\end{corollary}}

\def\geq{\geqslant}
\def\leq{\leqslant}
\def\ge{\geqslant}
\def\le{\leqslant}

 \def\R{\mathbb R}
 \def\R{\mathbb R}    
\def\N{\mathbb N}  
   
\def\<{\langle} \def\>{\rangle}

\def\0{{\mathbf{0}}}

\def\b{\beta}

\def\DB{{\mathfrak D}}

\allowdisplaybreaks

\title{Approximation and Interpolation in Kolmogorov-type groups}
\author{Antonello Pesce\thanks{Dipartimento di Matematica, Universit\`a di Bologna, Piazza di Porta s. Donato 5, 40126 Bologna, Italy. \textbf{e-mail}:
antonello.pesce2@unibo.it} }

\begin{document}
\maketitle

\begin{abstract}
We prove a real interpolation characterization for some non Euclidean H\"older spaces,  
built on the Lie structure induced by a class of ultra-parabolic Kolmogorov-type operators satisfying the H\"ormander condition.
As a by-product we also obtain an approximation property for intrinsically regular functions on the whole space.
\end{abstract}

\medskip


\section{Introduction}
\subsection{Statement of the problem}
We consider the problems of approximation and interpolation for some non Euclidean H\"older spaces built on the intrinsic geometry of the following Kolmogorov operator
\begin{equation}\label{PDE}
\mathcal{K}=\frac{1}{2}\sum_{i=1}^{p_0}\p_{x_i x_i}+Y, \quad Y=\langle Bx,\nabla\rangle +\p_t, \qquad (t,x)\in \R\times\R^{d},
\end{equation}
where $1\le p_0\le d$ and $B$ is a $d\times d$ constant matrix. If $p_0=d$ then $\mathcal{K}$ is a parabolic operator while in general, for $p_0<d$, $\mathcal{K}$ is degenerate, and not uniformly parabolic. However, $\mathcal{K}$ can be a hypoelliptic operator, provided some structural assumptions on the matrix $B$ are verified. This kind of PDEs appears in many linear and non-linear models in physics (see, for instance, \cite{Cercignani}, \cite{Desvillettes}) and in mathematical finance (see, for instance, \cite{BarucciPolidoroVespri}, \cite{CIBELLI201987}).

In \cite{LanconelliPolidoro}, Lanconelli and Polidoro first studied the non-Euclidean intrinsic Lie group structure induced by the vector fields $\p_{x_1},\dots, \p_{x_{p_0}},Y$, and from there on, many authors carried the study of general Kolmogorov operators in the framework of the theory of homogeneous groups. We recall, among many, \cite{DiFrancescoPascucci2}, who proved the existence of a fundamental solution under optimal regularity assumptions, \cite{DiFrancescoPolidoro}, \cite{Imbert} and \cite{MR1751429}, who proved Schauder type estimates for the associated Cauchy and Dirichlet problems, \cite{Imbert2}, who considered the case of integro-differential diffusions. 
One basic tool in the study of regularity properties of solutions to Kolmogorov operators with variable coefficients is the notion of H\"older spaces which are \textit{intrinsic} to the aforementioned geometric structure. Over the years, different notions of H\"older classes have been proposed by several authors, based on the Carnot-Carath\'eodory distance associated to the group (see, among many, \cite{DiFrancescoPolidoro}, \cite{MR1751429}, \cite{FrentzNystromPascucci2010}). Eventually Pagliarani et al. in \cite{MR3429628} established an optimal definition of $C_B^{n,\a}$ for any $n\in\N_0$, which is given in terms of the regularity in the directions of the vector fields $\p_{x_1}, \dots , \p_{x_{p_0}}$ and $Y$ with their different formal degrees: for instance,  the case $n=0$ consists of functions that are $C^\a$ continuous in the directions $\p_{x_1}, \dots , \p_{x_{p_0}}$ and $C^{\frac{\a}{2}}$ in $Y$. With this definition they managed to prove a Taylor-type formula with related estimates of the remainder, given in terms of the intrinsic distance.

Our main result, Theorem \ref{TH}, shows that these intrinsic spaces interpolate well with each other through the so called \textit{real method}, namely that
\begin{equation}\label{intro_interpolation}
\left(C_B^{n_1,\a_1},C_B^{n_2,\a_2}\right)_{\theta,\infty}=C_B^{n,\a},
\end{equation}
for any $n_1,n_2\in \N_0$, $\a_1,\a_2\in [0,1]$, for some expected parameters $n\in\N_0$, $\a\in (0,1)$.
In the Euclidean case, this identification is well known (see  \cite{MR3753604}, \cite{Triebel}), and can be proved by means of the general theory of semigroups. 
Similar arguments can be adapted for a different class of H\"older spaces also associated to degenerate elliptic and parabolic differential equations, the so called \textit{anisotropic} H\"older spaces, which, unlike the intrinsic spaces, only refer to the spatial variables (see Section \ref{Lett} for a more in-depth discussion). 
Here we choose to pursue a constructive approach based on explicit approximations of intrinsically regular functions (Theorem \ref{PROP}), which can possibly be of independent interest.
This approach is more or less inspired by the approximation theory for Sobolev functions (see \cite{Adams1975}, Chapter 5), 
it does not rely on the theory of semigroups, but instead, it heavily exploits the aforementioned intrinsic Taylor formula. 
\medskip

\noindent {Classically, characterizations of type \eqref{intro_interpolation} have remarkable applications in the study of optimal regularity 
in H\"older classes for linear elliptic and parabolic differential operators (see \cite{MR1406091}, \cite{MR3753604}, \cite{Lunardi2}, \cite{Lunardi}) as well as approximation theory (see \cite{MR0230022}). 
In the context of Kolmogorov operators, many works already exploit some kind of interpolation inequalities to obtain different versions of Schauder estimates in different settings (see, for instance \cite{Imbert}, \cite{DiFrancescoPolidoro} and \cite{Manfredini} among others, see also \cite{Lucertini2} for the most recent global estimate for the evolution problem, as well as a complete overview of the existing literature). 
Let us also mention that our results, together with the embeddings of intrinsic Sobolev spaces of integer order in \cite{PP22} could be applied to deduce Morrey-type embeddings for suitable notions of Sobolev-Slobodeckij and Besov intrinsic spaces as well, which, as far as our understanding, are yet to be explored in the literature in generality. Such studies will be subjects of future research. For the application of interpolation results to fractional embeddings in the Euclidean setting, see for instance, the references \cite{Adams1975} or \cite{MR2328004}.}
\medskip

The paper is organized as follows: in the remaining part of this Section, we specify the precise functional framework, state the structural hypothesis on the matrix $B$ and our main interpolation results, Proposition \ref{PROP2} and Theorem \ref{TH}, and finally discuss the previous literature on the subject and justify our approach.
In Section \ref{Preliminaries} we collect some results that are preliminary to the subsequent proofs; in Section \ref{approximation} we prove the core result of this work, the approximation property for intrinsically regular functions, which eventually leads to the proof of  our main results in Section \ref{proof}.

\subsection{H\"older spaces and real interpolation}\label{risultati}
We start by recalling the Lie group structure that is naturally associated to the differential operator \eqref{PDE}. In \cite{LanconelliPolidoro} the authors observed that $\mathcal{K}$ is invariant with respect to the left translations of the group $(\R^{d+1},\circ)$, where the non-commutative group law is defined by 
\begin{equation}\label{group_law}
(t,x)\circ (s,\x)=\ell_{(t,x)}(s,\x):=(s+t,\x+e^{sB}x), \qquad (t,x),\ (s,\x)\in\R^{d+1}.
\end{equation}  
Precisely, $\ell_{(t,x)}\circ\mathcal{K}=\mathcal{K}\circ\ell_{(t,x)}$ for any $(t,x)\in\R^{d+1}$.
Notice that $(\R^{d+1},\circ)$ is a group with identity $\text{Id}=(0,0)$ and inverse $(t,x)^{-1}=(-t,e^{-tB}x)$.

The authors also proved that the operator $\mathcal{K}$ is hypoelliptic if and only if, on a certain basis of $\R^{d+1}$ the matrix $B$ takes the form 
\begin{equation}\label{B_0}
B=\begin{pmatrix} \ast & \ast &\cdots & \ast & \ast \\
{\bf B}_1 &\ast &\cdots & \ast & \ast\\
0_{p_2\times p_0}& {\bf B}_2 &\cdots & \ast & \ast\\
\vdots &\vdots & \ddots &\vdots &\vdots \\
0_{p_r\times p_0}& 0_{p_r\times p_1} &\cdots & {\bf B}_r  & \ast
\end{pmatrix},
\end{equation}
where $\mathbf{B}_j$ is a $p_j\times p_{j-1}$ matrix of rank $p_j$, with $$p_0\geq p_1 \geq \dots \geq p_r\geq 1,\quad \sum_{j=1}^rp_j=d,$$
$0_{p_i\times p_j}$ is a $p_i\times p_j$ null matrix and the $\ast$-blocks are arbitrary. Throughout the paper we assume the following 
\begin{Assumption}\label{Ass}$B$ is a constant $d\times d$ matrix as in \eqref{B_0}, where each $\ast$-block is null.
\end{Assumption}
\noindent If (and only if) Assumption \ref{Ass} holds, then $\mathcal{K}$ is homogeneus of degree two with respect to the family of dilations $(D_{\l})_{\l>0}$ on $\R^{d+1}$ given by $$D_{\l}=(\l^2,\l I_{p_0}, \l^3 I_{p_1}\dots, \l^{2r+1}I_{p_r}),$$
where $I_{p_j}$ are $p_j\times p_j$ identity matrices, i.e. $(\mathcal{K} D_{\l}u)(t,x)=\l^{2}(\mathcal{K} u)(D_{\l}(t,x))$, for any $(t,x)\in\R^{d+1}$, $\l>0$. 
In this case, the matrix $B$ uniquely identifies the \textit{homogeneous Lie group} $\cG_B:=(\R^{d+1},\circ,D(\l))$.
We define the $D(\l)$-homogeneous quasi-norm on $\cG_B$ as follows: let
\begin{equation}\label{bar_p}
\bar{p}_j=p_0+\dots+p_j,\qquad j=0,\dots , r,
\end{equation}
and $\bar{p}_{-1}\equiv 0$. Then
\begin{equation}\label{norm}
\|(t,x)\|_B:=|t|^{\frac 12}+|x|_B, \qquad |x|_B:=\sum_{j=0}^r\sum_{i=\bar{p}_{j-1}+1}^{\bar{p}_j}|x_i|^{\frac{1}{2j+1}}.
\end{equation}
\smallskip

Next we recall the notions of $B$-intrinsic regularity and $B$-intrinsic H\"older space introduced in \cite{MR3429628}.
We start with some useful notations:
for any $(t,x)\in \R^{d+1}$, $i=1,\dots , p_0$ we denote 
\begin{equation}\label{vfields}
e^{\d\p_{x_i}}(t,x)=(t,x+\d e_i), \quad {e^{\d Y}(t,x)=(t+\d,e^{\d B}x)}, \qquad \d\in \R,
\end{equation}
the (unique) integral curves of the vector fields $\p_{x_i}$ and $Y$ respectively.
\begin{definition}
Let $X\in \{\p_{x_1},\dots, \p_{x_{p_0}},Y\}$ and $\a \in ]0,1]$. We say that $u\in C^{\a}_X$ if 
\begin{equation}
[u]_{C^\a_X}:=\sup_{z\in {\R^{d+1}}\atop \delta\in \R\setminus \{0\}}\frac{\left|u(e^{\delta X}z)-u(z)\right|}{|\delta|^{{\a}}}<\infty.
\end{equation}
We say that $u$ is $X$-differentiable in $z$ if the function $\d\mapsto u(e^{\d X}z)$ is differentiable in $0$ and refer to the function 
$z\mapsto Xu(z):=\frac{d}{d\d}u(e^{\d X}z)\mid_{\d=0}$ as the $X$\textit{-Lie derivative of u}.
\end{definition}

\begin{definition}\label{B-Holder}
Let $\a\in {[0,1]}$, then: 
\begin{itemize}
\item[i)] If $\a=0$, $C^{0,0}_B\equiv C(\R^{d+1})$ is the space of bounded and continuous function on $\R^{d+1}$, endowed with the sup-norm $$\|u\|_{C^{0,0}_B}:=|u|_{\infty}.$$ 
If $\a\in ]0,1]$ we say that $u\in C^{0,\a}_B$ if $u\in C^{\frac{\a}{2}}_Y$ and $u\in C^{\a}_{\p_{x_i}}$ for any $i=1,\dots , p_0$, with norm
\begin{align*}
\|u\|_{C^{0,\a}_B}&:=|u|_{\infty}+[u]_{C^{\frac{\alpha}{2}}_Y}+\sum_{i=1}^{p_0}[u]_{C^{\alpha}_{\p_{x_i}}}<\infty, \quad \a\in ]0,1].
\end{align*}
\item[ii)] 
$u\in C^{1,\alpha}_B$ if $u\in C^{\frac{1+\a}{2}}_Y$ and $\p_{x_i}u\in C^{\a}_B$ for any $i=1,\dots , p_0$, with norm
$$\|u\|_{C^{1,\alpha}_B}:=|u|_{\infty}+[u]_{C^{\frac{1+\alpha}{2}}_Y}+\sum_{i=1}^{p_0}\|\p_{x_i}u\|_{C^{0,\alpha}_B}<\infty.$$
\item[iii)] for $n\in\N$ with $n\geq 2$, $u\in C^{n,\alpha}_B$ if $Yu\in C^{n-2,\a}_B$ and $\p_{x_i}u\in C^{n-1,\a}_B$ for any $i=1,\dots , p_0$, with norm
$$\|u\|_{C^n_B}:=|u|_{\infty}+\|Yu\|_{C^{n-2,\alpha}_B}+\sum_{i=1}^{p_0}\|\p_{x_i}u\|_{C^{n-1,\alpha}_B}<\infty.$$
\end{itemize}
\end{definition}
\begin{remark}$C^{n,\a}_B$ is a Banach space. Notice also that $C^{n,\a}_{B}\subseteq C^{n',\a'}_{B}$ for $0\le n\le n'$ and $0\le \a'\le \a\le 1$.
\end{remark}
\medskip
\noindent For any two real Banach spaces $Z_1, Z_2$, by $Z_1=Z_2$ we mean that $Z_1$ and $Z_2$ have the same elements with equivalent norms and by $Z_1\subseteq Z_2$ we mean that $Z_1$ is continuously embedded in $Z_2$. We recall that the couple $\{Z_1, Z_2\}$ is called \textit{interpolation couple} if both $Z_1$ and $Z_2$ are continuously embedded in some Hausdorff topological vector space $\mathcal{Z}$; in this case, the intersection $Z_1\cap Z_2$ and the sum $Z_1+Z_2=\{u_1+u_2,\ u_1\in Z_1, \ u_2\in Z_2\}$, endowed with the norms 
$$\|u\|_{Z_1\cap Z_2}:=\max\{\|u\|_{Z_1},\|u\|_{Z_1}\}, \quad \|u\|_{Z_1+ Z_2}:=\inf\{\|u_1\|_{Z_1}+\|u_2\|_{Z_1}, u_1\in Z_1, \ u_2\in Z_2\},$$
are Banach spaces. Moreover, any Banach space $E$ such that  
\begin{equation}\label{intermediate}
Z_1\cap Z_2\subseteq E\subseteq Z_1+Z_2,
\end{equation}
is called an \textit{intermediate space}. Among these, we have the 
\textit{real interpolation space} $(Z_1,Z_2)_{\theta,\infty}$, which can be defined, for any  $\theta\in (0,1)$ (see, for instance \cite{MR3753604}), as 
\begin{equation}\label{k_interpolation}
(Z_1,Z_2)_{\theta,\infty}:=\{u\in Z_1+Z_2, \ \|u\|_{(Z_1,Z_2)_{\theta,\infty}}:=\sup_{\l\geq 0}\l^{-\theta}K(\l,u;Z_1,Z_2)<\infty\} 
\end{equation}
where 
$$K(\l,u;Z_1,Z_2):=\inf\left\{\|u_1\|_{Z_1}+\l\|u_2\|_{Z_2},\, u=u_1+u_2, \, u_1\in Z_1,\, u_2\in Z_2\right\}.$$  
\medskip

\noindent Our first results states that any space $C^{n,0}_B$ belongs to a special class of \textit{intermediate spaces} between the larger space $C^{n_1,0}_B$, $n_1<n$ and a smaller space $C^{n_2,0}_B$, $n<n_2$, with  $n,n_1,n_2\in\N_0$. 
\begin{proposition}\label{PROP2}
Let $n,n_1,n_2\in\N_0$, $n_1<n<n_2$.
\begin{itemize}
\item[i)]There exists a constant $c_1$ such that, for any $u\in C^{n_2,0}_B$, we have 
\begin{equation}
\|u\|_{C^{n,0}_B}\le c_1 \|u\|^{\frac{n_2-n}{n_2-n_1}}_{C^{n_1,0}_B}\|u\|^{\frac{n-n_1}{n_2-n_1}}_{C^{n_2,0}_B}; 
\end{equation}
\item[i)]There exists a constant $c_2$ such that, for any $u\in C^{n,0}_B$, we have 
\begin{equation}
K(\l,u;C^{n_1,0}_B,C^{n_2,0}_B)\le c_2 \l^{\frac{n-n_1}{n_2-n_1}}\|u\|_{C^{n,0}_B}, \quad \l>0.
\end{equation}
\end{itemize}
Equivalently (see, for instance \cite{MR3753604}, Proposition 1.20) it holds that 
 \begin{equation}
\left(C^{n_1,0}_B,C^{n_2,0}_B\right)_{\frac{n-n_1}{n_2-n_1},1}\subseteq C^{n,0}_B\subseteq \left(C^{n_1,0}_B,C^{n_2,0}_B\right)_{\frac{n-n_1}{n_2-n_1},\infty}.
\end{equation}
\end{proposition} 
\noindent Next, our main result, gives a characterization of $B$-H\"older spaces of any order as interpolation spaces.
\begin{theorem}\label{TH}
Let $n_1,n_2\in\N_0$ and $\a_1,\a_2\in [0,1]$, with $n_1+\a_1\le n_2+\a_2$, and let $\theta\in (0,1)$. Let $n\in\N_0$, $n_1\le n\le n_2$ such that
\begin{equation}\label{eTHbis}
\frac{n-(n_1+\a_1)}{(n_2+\a_2)-(n_1+\a_1)}<\theta <\frac{n+1-(n_1+\a_1)}{(n_2+\a_2)-(n_1+\a_1)}.
\end{equation}
If $\a:=(n_1+\a_1)+\theta[n_2+\a_2-(n_1+\a_1)]-n \neq 0,1$
then it holds that
\begin{equation}\label{eTH}
\left(C^{n_1,\a_1}_B,C^{n_2,\a_2}_B\right)_{\theta,\infty}=C^{n,\alpha}_B
\end{equation}
\end{theorem}
\begin{remark}\label{open}If $\Omega$ is an open set of $\R^{d+1}$, one can define the H\"older spaces $C^{n,\a}_B(\Omega)$, $n\in\N_0$, $\a\in (0,1)$ as in the case $\Omega=\R^{d+1}$, where in particular, for $X\in\{\p_{x_1},\dots \p_{x_{p_0}},Y\}$,
\begin{equation*}
[u]_{C^{\a}_X(\Omega)}:=\sup\left\{\frac{|u(e^{\d X}z)-u(z)|}{|\d|^{\a}}, \ z\in\Omega, \ \d\in \R\setminus\{0\}\text{ s.t. } e^{\d'X}z\in\Omega,\ |\d'|<|\d|\right\}.
\end{equation*}
Then it is {likely} that our results also hold for the spaces $C^{n,\a}_B(\Omega)$, provided that the domain $\Omega$ is smooth enough (see \cite{Triebel}, Theorem 1.2.4 and Section 4.5). The construction of the appropriate extension operators under optimal regularity assumptions for $\Omega$ is outside the scope of the present paper and is postponed to future research.
\end{remark}
\begin{remark}
Notice that the assumption $n_1+\a_1\le n_2+\a_2$ is not restrictive. Indeed we can always write 
$$\left(C^{n_1,\a_1}_B,C^{n_2,\a_2}_B\right)_{\theta,\infty}=\left(C^{n_2,\a_2}_B,C^{n_1,\a_1}_B\right)_{1-\theta,\infty}.$$
\end{remark}

\subsection{Related works and strategy of the proof}\label{Lett}
As already mentioned, in the Euclidean case, the interpolation identity \eqref{intro_interpolation} can be proved by means of the theory of semigroups (\cite{MR3753604}, \cite{Triebel}): 
let $C, C^1$ denote the spaces of bounded continuous functions and bounded continuously differentiable functions with bounded derivatives respectively, and let $C^{\theta}$ denote the space of bounded and uniformly $\theta$-H\"older continuous
functions, all endowed with the usual norms.  Roughly speaking, denoting with $T_i$ the translation semigroup of the $i$-th variable, and $\cD(\p_{x_i})=\{f\in C, \ \exists \p_{x_i}f\in C\}$, we have 
\begin{equation}\label{semi_approach}
(C,C^1)_{\theta,p}=\Big(C,\bigcap_{i=1}^d\cD(\p_{x_i})\Big)_{\theta,p}=\bigcap_{i=1}^d\Big(C,\cD(\p_{x_i})\Big)_{\theta,p},
\end{equation}
for any $p\ge 1$, where 
\begin{equation}\label{semi2}
\big(C,\cD(\p_{x_i})\big)_{\theta,\infty}=\{f\in C, \ t\mapsto \phi(t):=t^{-\theta}\|T_i(t)f-f\|_{\infty}\in L^{\infty}(0,\infty)\},
\end{equation}
that is the space of $\theta$-H\"older continuous functions in the $i$-th variable. Therefore $(C,C^1)_{\theta,\infty}=C^{\theta}$. 
In particular \eqref{semi2} holds more generally for $p\ge 1$, and substituting $T_i$ and $\p_{x_i}$ with any bounded semigroup $T$ and its infinitesimal generator on a Banach space, while the second equality in \eqref{semi_approach} precisely holds because the semigroups $T_i$, $i=1,\dots, d$ commute. 

In \cite{Lunardi}, Lunardi adapted these arguments to derive an interpolation result for anisotropic H\"older spaces associated with operators of the type \eqref{PDE}, which are spaces of functions on $\R^d$ (only the spatial variables are considered) defined with respect to the metric induced by $|\cdot|_B$ (cfr. \eqref{norm}). For instance, let us consider the degenerate prototype of \eqref{PDE}, that is the Langevin operator 
$$\mathcal{K}=\frac 12\p_{x_1x_1}+Y, \quad Y=x_1\p_{x_2}+\p_t,$$
which corresponds to the matrix 
\begin{equation*}
B=\begin{pmatrix} 0 & 0 \\ 1 &0 \end{pmatrix}.
\end{equation*}
In this case the anisotropic H\"older regularity in the space variables reads: 
\begin{equation}\label{anisotropic}
|u(x)- u(\x)|\le C|x-\x|^\a_B\sim |x_1-\x_1|^{\a}+|x_2-\x_2|^{\frac{\a}{3}}, \quad x,\x\in\R^2.
\end{equation}
A regularity theory for equations with variable coefficients, built on this functional setting, is useful when one needs to consider coefficients or data that are irregular in time (see also \cite{MR4358660}, \cite{MR4355925}, \cite{CdRMPZ22}, \cite{Lucertini} for recent developments). 
This approach leads to mild or weak/distributional solutions also because it 
cannot benefit from any regularizing effect of the semigroup in time. 

On the other hand, for $B$ as above, a function $u\in C^{0,\a}_B$ satisfies (see Theorem \ref{Taylor} in the next Section), for any $(t,x),(s,\x)\in\R^3$, $$|u(t,x)- u(s,\x)|\le C\|(s,\x)^{-1}\circ (t,x)\|^\a_B\sim |t-s|^{\frac{\a}{2}}+|x_1-\x_1|^{\a}+|x_2-\x_2-(t-s)\x_1|^{\frac{\a}{3}}.$$
Notice that for points that are on the same time level, the increments are controlled as for \eqref{anisotropic}; otherwise, as opposed to the standard parabolic case, the regularity in the euclidean directions is somehow entangled due to the group law associated to $B$.  
This functional setting, which exploits some regularity in the full set of variables, is more suited to obtain classical solutions (see \cite{DiFrancescoPolidoro}, \cite{MR1751429} and the references therein). To the best of our knowledge this is the first attempt at obtaining a general interpolation characterization for the spaces $C^{n,\a}_B$. For this purpose, notice that while it is true that $$\big(C,\cD(Y)\big)_{\theta,\infty}=\{f\in C, \ t\mapsto \phi(t):=t^{-\theta}\|e^{tY}f-f\|_{\infty}\in L^{\infty}(0,\infty)\}=C_Y^{\theta},$$ since the vector fields $\p_{x_1}$ and $Y$ don't commute, it's not at all clear whether
\begin{equation}\label{question}
\big(C,\cD(\p^2_{x_1})\cap \cD(Y)\big)_{\frac{\theta}{2},\infty}=\left(C,\cD(\p^2_{x_1})\right)_{\frac{\theta}{2},\infty}\cap \left(C, \cD(Y)\right)_{\frac{\theta}{2},\infty}. 
\end{equation}
This makes it difficult to pursue an approach based on semigroup theory.
Instead of trying to prove \eqref{question} we follow a constructive approach based on approximation: for instance, for any $u\in C^{0,\a}_B$ we can define $$u_{\eps}(z)=\int_{\R^{d+1}}u(\z)\phi\left(D_{\eps^{-1}}(\z^{-1}\circ z)\right)\frac{d\z}{\eps^{Q+2}}, \quad z\in \R^{d+1},$$
where $\phi$ is a test function with unit integral, supported on $\|y\|_B\le 1$. This approximation complies with the Lie structure induced by $B$, and indeed  we prove the following expected controls
\begin{equation}\label{stime}
|u-u_{\eps}|_{\infty}\le c\eps^{\a}\|u\|_{C^{0,\a}_B}, \quad \|u_{\eps}\|_{C^{1,0}_B}\le c\eps^{\a-1}\|u\|_{C^{0,\a}_B},
\end{equation}
so that 
$$\eps^{-\a}K(\eps,u;C,C^{1,0}_B)\le \eps^{-\a}(|u-u_{\eps}|_{\infty}+\eps\|u_{\eps}\|_{C^{1,0}_B})\le 2c \|u\|_{C^{0,\a}_B}.$$
This way of reasoning is fairly standard (see, for instance, Chapter 1 in \cite{MR3753604}), though the estimates \eqref{stime} are new for the intrinsic spaces and we present a generalization for any order $n\in\N_0$ in Theorem \ref{PROP}.

\section{Preliminaries}\label{Preliminaries}
Here we collect some remarks and known preliminary results which will play a central role for the next Sections. 
\begin{remark}\label{r1}
Under Assumption \ref{Ass} the matrix $B$ is a nilpotent matrix of degree $r+1$. In particular 
 $$e^{\delta B}=I_d+\sum_{j=1}^r\frac{B^j}{j!}\delta^j,$$ 
is a lower triangular matrix with diagonal $(1,\dots, 1)$ and therefore it has determinant equal to 1. 
Indeed, by Assumption \ref{Ass}, a direct computation shows that, for any $n\leq r$ (recalling the notation \eqref{bar_p} for $\bar{p}_j$) we have 
\begin{equation}\label{Bn}
B^n=\begin{pmatrix} 0_{\bar{p}_{n-1}\times p_0}& 0_{\bar{p}_{n-1}\times p_1} &\cdots & 0_{\bar{p}_{n-1}\times p_{r-n}}& 0_{\bar{p}_{n-1}\times (\bar{p}_{r}-\bar{p}_{r-n}})\\
\prod_{j=1}^n{\bf B}_j &0_{p_n\times p_1} &\cdots & 0_{p_n\times p_{r-n}}& 0_{p_n\times (\bar{p}_{r}-\bar{p}_{r-n})}\\
0_{p_{n+1}\times p_0}& \prod_{j=2}^{n+1}{\bf B}_j &\cdots & 0_{p_{n+1}\times p_{r-n}}& 0_{p_{n+1}\times (\bar{p}_{r}-\bar{p}_{r-n})}\\
\vdots &\vdots & \ddots &\vdots &\vdots \\
0_{p_r\times p_0}& 0_{p_r\times p_1} &\cdots & \prod_{j=r-n+1}^{r}{\bf B}_j  & 0_{p_r\times (\bar{p}_{r}-\bar{p}_{r-n})}
\end{pmatrix}, \qquad \prod_{j=1}^{n}{\bf B}_j={\bf B}_n{\bf B}_{n-1}\cdots {\bf B}_1,
\end{equation}
and ${B}^n=0$ for $n>r$. 
\end{remark}

\begin{remark}\label{r3}
A simple computation shows that, for any $z,y\in \R^{d+1}$,
\begin{equation}\label{invariance_law}
y^{-1}\circ e^{\delta Y}z=e^{\delta Y}(y^{-1}\circ z), \qquad y^{-1}\circ e^{\delta \p_{x_i}}z=e^{\delta \p_{x_i}}(y^{-1}\circ z), 
\quad i=1,\dots , d,
\end{equation}
which is a way to restate the invariance of operator \eqref{PDE} w.r.t the group law.
\end{remark}

\begin{lemma}\label{r2}
The following identities hold:
\begin{equation}\label{homogeneity}
D_\l e^{\delta Y}(z)=e^{\delta\lambda^2Y}\left(D_\l z\right), \qquad D_\l e^{\delta \p_{x_i}}(z)=e^{\delta \l^{2j+1} \p_{x_i}}D_\l z, 
\quad i=\bar{p}_{j-1}+1,\dots , \bar{p}_{j},
\end{equation}
for any $z\in\R^{d+1}$ and $\d,\l\in\R$. 
\end{lemma}
\proof Let $z=(t,x)$, and let $D^0_\l$ be the spatial part of the dilation operator $D_\l$, namely $D^0_{\l}=(\l I_{p_0}, \l^3 I_{p_1},\dots, \l^{2r+1}I_{p_r})$. Denoting with $x_{\{i\}}$, $i=0,\dots , r$, the projection of $x$ on $\R^{p_i}$ defined by 
$$x_{\{i\}}:=(x_{\bar{p}_{i-1}+1},\dots , x_{\bar{p}_{i}}),$$
by \eqref{Bn} we can write 
\begin{equation}\label{Bflow}
\left(e^{\delta B}x \right)_{\{i\}}=x_{\{i\}}+\sum_{j=0}^{i-1}\delta^{i-j}\left(\prod_{k=j+1}^{i}{\bf B}_k\right)x_{\{j\}}, \qquad i=0,\dots , r.
\end{equation}
Then we have
\begin{align*}
\left(D^0_{\l}e^{\delta B}x\right)_{\{i\}}&=\l^{2i-1}x_{\{i\}}+\l^{2i-1}\sum_{j=0}^{i-1}\delta^{i-j}\left(\prod_{k=j+1}^{i}{\bf B}_k\right)x_{\{j\}}\\
&=\l^{2i-1}x_{\{i\}}+\sum_{j=0}^{i-1} (\l^2\delta)^{i-j}\left(\prod_{k=j+1}^{i}{\bf B}_k\right)\l^{2j-1}x_{\{j\}}
=\left(e^{\delta\lambda^2B}\left(D^0_{\l}x\right)\right)_{\{i\}}.
\end{align*}
where in the last equality we used that $\l^{2i-1}x_{\{i\}}=(D_\l x)_{\{i\}}$ and $\l^{2j-1}x_{\{j\}}=(D_\l x)_{\{j\}}$. Thus we get
$$D_{\l}e^{\delta Y}z=\left(\lambda^2(t+\delta),e^{\delta\lambda^2B}\left(D^0_{\l}x\right)\right)=e^{\delta\lambda^2Y}\left(D(\lambda)z\right).$$
The second equality is clear from the definitions.
\endproof

We conclude this Section by recalling the crucial intrinsic Taylor expansion proved in \cite{MR3429628}. Here $\b=(\b_1,\dots, \b_d)\in \N^d_0$ denotes a multi-index, and as usual $\b !:=\prod_{j=1}^d(\b_j!)$. Moreover, for any $x\in\R^d$, $x^\b=x_1^{\b_1}\dots , x_d^{\b_d}$ and $\p^{\b}=\p_x^{\b}=\p_{x_1}^{\b_1}\cdots \p_{x_d}^{\b_d}$. The \textit{$B$-length} of $\b$ is defined as 
$$|\b|_B:=\sum_{i=0}^r(2i+1)\sum_{j=\bar{p}_{i-1}+1}^{\bar{p}_i}\b_j.$$ 

\begin{theorem}[Intrinsic Taylor formula]\label{Taylor} If $u\in C^{n,\a}_B$ then we have 
\begin{equation}\label{B-regularity}
Y^k\p^\b_xu\in C_B^{n-2k-|\b|_B,\a}, \qquad 0\le 2k+|\b|_B\le n,
\end{equation}
and 
\begin{equation}\label{T_remainder}
|u(z)-T_nu(\z,z)|\le c_B\|u\|_{C_B^{n,\a}}\|\z^{-1}\circ z\|^{n+\a}_B, \qquad z,\z\in\R^{d+1},
\end{equation}
where $c_B$ is a positive constant which only depend on $B$ and $T_nu(\z,\cdot)$ is the $n$-th order \textit{$B$-Taylor polynomial} of $u$ around $\z=(s,\x)$ defined as
\begin{equation}\label{Polynomial}
T_nu(\z,z)=\sum_{0\le 2k+|\b|_B\le n}\frac{1}{k!\b!}\left(Y^k\p^\b_\x u(s,\x)\right)(t-s)^k(x-e^{(t-s)B}\x)^\b, \qquad z=(t,x)\in \R^{d+1}.
\end{equation}
\end{theorem}
\noindent\textbf{Notation}: Eventually, we will frequently use the notation $\lesssim$. For two quantities $Q_1$ and $Q_2$, we mean by $Q_1 \lesssim Q_2$ that there exists a constant $c_B$, only dependent on the matrix $B$, such that $Q_1\le c_B Q_2$.

\section{Approximation of $C^{n,\a}_B$ functions}\label{approximation}
This Section is mainly devoted to the proof of the following Theorem \ref{PROP}, which establishes an approximation property for our $B$-H\"older continuous functions.

\begin{theorem}[Approximation property]\label{PROP}
Let $n,l \in \N_0$ and $\a\in [0,1]$. There exists a constant $c_B$, only dependent on $B$ and $d$, such that for any $u\in C^{n,\a}_B$ and $0<\eps \le 1$, there exists $u_{\eps}\in C^{\infty}$ such that
\begin{align}\label{Convergence}
\|u-u_{\eps}\|_{C^{l,0}_B}&\leq c_B \eps^{n+\a-l}\|u\|_{C^{n,\a}_B}, & l\le n,\\
\|u_{\eps}\|_{C^{l,0}_B}&\leq 
c_B \eps^{n+\a-l}\|u\|_{C^{n,\a}_B} \quad  & {l>n}.
\end{align}
\end{theorem}


The proof is based on the intrinsic H\"older expansion of Theorem \ref{Taylor} and the invariance properties of Lemma \ref{r2} and Remark \ref{r3}. We also need the following technical Lemma, which remarkably says that we can exchange the differentiation of the intrinsic Taylor polynomial with the differentiation of the inner function, taking into account the intrinsic degree of the derivatives $\p_{i}$, $i=1,\dots , p_0$ and $Y$. 

To avoid ambiguity, when necessary, we will write $Y_z$, $Y_\z$ as well as $\p_{x_i}$, $\p_{\x_i}$ to clearly distinguish the variables with respect to which the derivation is applied. 

\begin{lemma}\label{LEM}
For any $u\in C^{n,\a}_B$ and $z=(t,x), \, \z=(s,\x)\in \R^{d+1}$ we have 
\begin{align}
\p_{x_i}T_nu(\z,z)&=T_{n-1}(\p_{i}u)(\z,z), \quad n\ge 1, \quad i=1,\dots , p_0, \label{Taylor1}\\
Y_zT_nu(\z,z)&=T_{n-2}(Yu)(\z,z), \quad n\ge 2. \label{Taylor2}
\end{align}
\end{lemma}
\proof Equality \eqref{Taylor1} is analogous to the Euclidean case and follows by a direct computation. Indeed it suffices to notice that, for $i=1,\dots, p_0$, $T_nu(\z,z)$ depends on $x_i$ only in the terms $(x-e^{(t-s)B}\z)_i$ and $\p_{x_j}(x-e^{(t-s)B}\z)_i=\d_{ij}$ ($\d_{ij}$ being the Kronecker Delta) for $i,j=1,\dots, p_0$. 
Equality \eqref{Taylor2} on the other hand, needs a more careful analysis: the directional derivative $Y_z$ is not null when applied both to the terms $(t-s)^k$ and $(x-e^{(t-s)B}\z)_i$ for any $i>p_0$ and, most importantly, some commutators need to be introduced in order to exchange the order of derivation between $\p_{\x_i}$ and $Y_\z$, and make $Y_\z u$ appear in any term of the sum.
\medskip

First we introduce some notations: let $v:=x-e^{(t-s)B}\x\in\R^{d}$ and let us denote with $v^{[i]}$, $i=0,\dots, r$, the projection of $v$ on $\{0\}^{\bar{p}_{i-1}}\times \R^{{p}_i}\times \{0\}^{d-\bar{p}_i}$; then we denote with $X_i$ the vector field $X_i:=\langle v^{[i]},\nabla\rangle$, were $\nabla$ represent the usual Nabla operator. Clearly we have
\begin{align}
Y_zX_0&=0, \label{eTaylor1}
\end{align}
Moreover we have
\begin{align}
Y_zX_iu(\z)&=\sum_{j=\bar{p}_{i-1}+1}^{\bar{p}_i}(Bv)_j\partial_{\x_j}u(\z)
=\langle (Bv)^{[i]},\nabla u(\z)\rangle \nonumber\\
&=\langle Bv^{[i-1]},\nabla u(\z)\rangle
=\left[X_{i-1},Y_{\z}\right]u(\z)=:X^{(1)}_{i-1}u(\z). \label{eTaylor2}
\end{align} 
\medskip
To prove \eqref{Taylor2} we proceed by induction. For $n=2$ we directly get
\begin{align*}
Y_zT_2u(\z,z)&=Y_z\left(u(\z)+X_0u(\z)+\frac{1}{2!}X_0^2u(\z)+Y_\z u(\z)(t-s)\right)
=Y_zY_\z u(\z)(t-s)=Y_\z u(\z).
\end{align*}
Then we proceed in two steps, namely, the induction from $2n+1$ to $2n+2$ and the induction from $2n$ to $2n+1$. 
\medskip\\
\textit{1: Induction from $2n+1$ to $2n+2$}
\medskip\\
Let $\b=(\b_0,\dots, \b_r)\in \N_0^{r+1}$ and let us denote just in the context of this proof $$|\b|_B=\sum_{i=0}^r(2i+1)\b_i, \quad\b!=\b_1!\cdots\b_r!.$$ Then, with the notation introduced above we can conveniently write 
\begin{align}
T_{2n+2}u(\z,z)&=T_{2n+1}u(\z,z)+\sum_{k=0}^{n+1}\sum_{|\b|_B=2(n+1)-2k}\frac{1}{k!\b!}Y_\x^k X_0^{\b_0}\cdots X_{n-k}^{\b_{n-k}}u(\z)(t-s)^k\\
&=:T_{2n+1}u(\z,z)+\tilde T_{2n+1}u(\z,z).
\end{align} 
By \eqref{eTaylor1} and \eqref{eTaylor2} we have 
\begin{align}
Y_zY_{\z}^k(t-s)^k&=Y_{\z}^k(t-s)^{k-1}, \label{eTaylor3}\\
Y_zX^{\b_i}_i&=\b_iX^{\b_i-1}_iX^{(1)}_{i-1} \label{eTaylor4}\\
Y_{\z}X^{\b_i}_i&=X^{\b_i}_iY_\z-\b_iX^{\b_i-1}_iX^{(1)}_{i}. \label{eTaylor5}
\end{align}
Therefore, applying $Y_z$ to $\tilde T_{2n+1}u(\z,z)$ we get $Y_z\tilde T_{2n+1}u(\z,z)=\sum_{k=0}^{n+1}S_k$ with 
\begin{align*}
S_0&=\sum_{|\beta|_B=2n+2}\frac{1}{\b!}\sum_{i=0}^n {\b_i} X_0^{\b_0}\cdots X^{\b_i-1}_{i}\cdots X^{\b_n}_nX^{(1)}_{i-1}u(\z)\\
&=\sum_{|\beta|_B=2n+2}\sum_{i=0}^{n-1}\frac{1}{\b_0!\cdots (b_{i+1}-1)!\cdots \b_n!}X_0^{\b_0}\cdots X^{\b_{i+1}-1}_{i+1}\cdots X^{\b_n}_nX^{(1)}_{i}u(\z)\\
&=\sum_{i=0}^{n-1}\sum_{|\beta|_B=2(n-i)-1}\frac{1}{\b!}X_0^{\b_0}\cdots X^{\b_{n-i-1}}_{n-i-1}X^{(1)}_{i}u(\z), 
\end{align*}
and, for $k=1,\dots , n+1$, 
\begin{align*}
 S_k=&\sum_{|\b|_B=2(n+1)-2k}\frac{(t-s)^{k-1}}{(k-1)!\b!} Y_{\z}^{k-1}\left( X_0^{\b_0}\cdots X_{n-k}^{\b_{n-k}}Y_\z-
 \sum_{i=0}^{n-k}\b_i X_0^{\b_0}\cdots X_{i}^{\b_{i}-1}\cdots X_{n-k}^{\b_{n-k}}X^{(1)}_i\right) u(\z)\\
 &\quad + \sum_{|\b|_B=2(n+1)-2k}\frac{(t-s)^k}{k!\b!}Y_\z^k\sum_{i=0}^{n-k}\b_iX_0^{\b_0}\cdots X^{\b_i-1}_{i}\cdots X^{\b_{n-k}}_{n-k}X^{(1)}_{i-1}u(\z)\\
 =&S_k^{(1)}-\sum_{i=0}^{n-k}\sum_{|\b|_B=2(n-k-i)+1}\frac{1}{(k-1)!\b!}Y_{\z}^{k-1}X_0^{\b_0}\cdots X^{\b_{n-i-k}}_{n-i-k}X^{(1)}_i u(\z) (t-s)^{k-1}\\
 &\quad + \sum_{i=0}^{n-k-1}\sum_{|\b|_B=2(n-k-i)-1}\frac{1}{k!\b!}Y_{\z}^{k}X_0^{\b_0}\cdots X^{\b_{n-i-k-1}}_{n-i-k-1 }X^{(1)}_i u(\z) (t-s)^{k}\\
 =&S_k^{(1)}+S_k^{(2)}+S_k^{(3)}.
\end{align*}
Now notice that
\begin{equation}\label{eTaylor6}
S_0+S_1^{(2)}=0\qquad S_k^{(3)}+S_{k+1}^{(2)}=0, \quad k=1,\dots, n. 
\end{equation}
Moreover $S_n^{(3)}=S_{n+1}^{(2)}=S_{n+1}^{(3)}=0$ by \eqref{eTaylor3}. Thus we finally get
\begin{align*}
Y_z\tilde{T}_{2n+2}u(\z,z)&=\sum_{k=0}^{n}S^{(1)}_{k+1}
=\sum_{k=0}^{n}\sum_{|\b|_B=2n-2k}\frac{1}{k!\b!}Y_{\z}^{k}X_0^{\b_0}\cdots X_{n-k-1}^{\b_{n-k-1}}Y_\z u(\z)(t-s)^k
=\tilde{T}_{2n}Yu(\z,z),
\end{align*}
and therefore, using the induction hypotheses
$$Y_zT_{2n+2}u(\z,z)=T_{2n-1}Yu(\z,z)+\tilde{T}_{2n}Yu(\z,z)=T_{2n}Yu(\z,z).$$
\medskip\\
\textit{2: Induction from $2n$ to $2n+1$}
\medskip\\
This Step only differs from the previous one by the tracking of the indexes throughout the computations and the fact that the polynomial $\tilde{T}_{2n+1}$ always contains at least one field of type $X_i$. We write 
$Y_z\tilde{T}_{2n+1}u(\z,z)=\sum_{k=0}^nS_k$ and we derive, with similar computations than in the previous step
$$ S_0=\sum_{i=0}^{n-1}\sum_{|\b|_B=2(n-i)}\frac{1}{\b!}X_0^{\b_0}\cdots X_{n-i-1}^{\b_{n-i-1}}X^{(1)}_iu(\z),$$
and, for any $k=1,\dots, n-1$, 
\begin{align*}
S_k=&\sum_{|\b|_{B}=2(n-k)+1}\frac{1}{(k-1)!\b!}Y_{\z}^{k-1}X_0^{\b_0}\cdots X_{n-k}^{\b_{n-k}}Y_\z u(\z)(t-s)^{k-1}\\
&\quad -\sum_{i=0}^{n-k}\sum_{|\b|_B=2(n-k-i)}\frac{1}{(k-1)!\b!}Y_{\z}^{k-1}X_0^{\b_0}\cdots X^{\b_{n-i-k}}_{n-i-k}X^{(1)}_i u(\z) (t-s)^{k-1}\\
&\quad + \sum_{i=0}^{n-k-1}\sum_{|\b|_B=2(n-k-i-1)}\frac{1}{k!\b!}Y_{\z}^{k}X_0^{\b_0}\cdots X^{\b_{n-i-k-1}}_{n-i-k-1 }X^{(1)}_i u(\z) (t-s)^{k}\\
 =&S_k^{(1)}+S_k^{(2)}+S_k^{(3)}.
\end{align*}
In particular, $S_0+S_1^{(2)}=0$, $S_k^{(3)}+S_{k+1}^{(2)}=0$ for any $k=1,\dots, n-1$ and $S_n^{(3)}=0$ by \eqref{eTaylor1}. Therefore we can eventually derive \eqref{Taylor2}.
\endproof
We are now ready to prove Theorem \ref{PROP}.
\proof[Proof of Theorem \ref{PROP}] 
Let $u\in C^{n,\a}_B$ and let $\phi$ be a test function supported on $\|z\|_B\le 1$ with unit integral. We define our candidate approximation by
\begin{align*}
u^{(n)}_{\eps}(z)=T_nu(\cdot,z)\star_B\phi_{\eps}(z):=\int_{\R^{d+1}}T_nu(\z,z)\phi\left(D_{\eps^{-1}}(\z^{-1}\circ z)\right)\frac{d\z}{\eps^{Q+2}}, 
\end{align*}
where $T_nu(\z,z)$ is the Taylor polynomial in \eqref{Polynomial} and $Q$ is the spatial \textit{homogeneous dimension} of $\R^{d+1}$ w.r.t $(D_\l)_{\l>0}$, that is the positive integer $$Q=p_0+3p_1+\cdots (2r+1)p_r.$$ Notice that, by the change of variables $\bar{\z}=D_{\eps^{-1}}(\z^{-1}\circ z)$, recalling Remark \ref{r1}, it is easy to verify that
\begin{equation}\label{mollifier}
\int_{\R^{d+1}}\phi\left(D_{\eps^{-1}}(\z^{-1}\circ z)\right)\frac{d\z}{\eps^{Q+2}}=\int_{\|\z\|_B\le 1}\phi(\z)d\z=1.
\end{equation}
In particular we have 
\begin{align*}
u(z)-u^{(n)}_{\eps}(z)=&\int_{\R^{d+1}}\left(u(z)-T_nu(\z,z)\right)\phi\left(D_{\eps^{-1}}(\z^{-1}\circ z)\right)\frac{d\z}{\eps^{Q+2}}.
\end{align*}
Observe now that the differentiation of $u^{(n)}_{\eps}$ falls back both on the test function and the polynomial $T_nu$, the latter being handled by Lemma \ref{LEM}; for the former case, by Lemma \ref{r2} and Remark \ref{r3} we deduce 
\begin{align*}
Y_z\phi (D_{\eps^{-1}}(\z^{-1}\circ z))&=\eps^{-2}Y\phi(\bar \z)\mid_{\bar \z=D_{\eps^{-1}}(\z^{-1}\circ z)},\\
\p_{z_j}\phi(D_{\eps^{-1}}(\z^{-1}\circ z))&=\eps^{-(2i+1)}\p_j\phi(\bar \z)\mid_{\bar\z=D_{\eps^{-1}}(\z^{-1}\circ z)}, \quad \bar p_{i-1}<j\leq \bar p_{i},
\end{align*}
and, more generally, for any $\beta\in \N_0^{d}$, 
\begin{equation}\label{eP1}
Y^k_z\p_z^{\beta}\phi (D_{\eps^{-1}}(\z^{-1}\circ z))=\eps^{-2k-|\beta|_B}Y^k\p^{\beta}\phi(\bar z)\mid_{\bar z=D_{\eps^{-1}}(\z^{-1}\circ z)}.
\end{equation}
\medskip
\textit{1: Preliminary controls}
\medskip\\
Let us denote a general intrinsic derivative of order $m$, $\DB^m=\sum_{2k+|\b|_B=m}Y^k\p^\b$, and let
$$I^{(n,m)}_{\eps}u(z):=\int_{\R^{d+1}}\left(u(z)-T_nu(\z,z)\right)\DB_z^{\beta}\phi (D_{\eps^{-1}}(\z^{-1}\circ z))\frac{d\z}{\eps^{Q+2}}.$$
We want to prove that
\begin{align}
|I^{(n,m)}_{\eps}u|_{\infty}&\lesssim \eps^{n+\a-m}\|u\|_{C^{n,\a}_B},\label{eP2}\\
[I^{(n,m)}_{\eps}u]_{C^{\frac{1}{2}}_Y}&\lesssim\eps^{n+\a-(m+1)}\|u\|_{C^{n,\a}_B}, \quad n\ge 1.\label{eP2bis}
\end{align}
By \eqref{eP1}, the change of variable $\bar\z=D_{\eps^{-1}}(\z^{-1}\circ z)$, and Theorem \ref{Taylor}, we have
\begin{align*}
|I^{(n,m)}_{\eps}u|_{\infty}&\le \int_{\R^{d+1}}\left|\left(u(z)-T_nu(\z,z)\right)\DB^m_z\phi (D_{\eps^{-1}}(\z^{-1}\circ z))\right|\frac{d\z}{\eps^{Q+2}}\\
&\lesssim \|u\|_{C^{n,\a}_B}\int_{\R^{d+1}}\|\z^{-1}\circ z\|_{B}^{n+\a}\eps^{-m}\DB^m\phi(\bar \z)\mid_{\bar \z=D_{\eps^{-1}}(\z^{-1}\circ z)}\frac{d\z}{\eps^{Q+2}}\\
&\lesssim \|u\|_{C^{n,\a}_B}\int_{\|\bar\z\|_B\le 1}\|D_{\eps}\bar \z\|_{B}^{n+\a}\eps^{-m}\DB^m\phi(\bar \z)d\bar \z\lesssim \eps^{n+\a-m}\|u\|_{C^{n,\a}_B},
\end{align*}
were we used that $\|D_{\eps}\bar \z\|_{B}=\eps\|\bar\z\|_B$ in the last inequality, and this proves \eqref{eP2}. On the other hand we have 
\begin{align*}
&|I^{(n,m)}_{\eps}u(z)-I^{(n,m)}_{\eps}u(e^{\d Y}z)|\\
&\qquad \le \int_{\R^{d+1}}\left|u(z)-T_nu(\z,z)\right|
\left|\DB_{z'}^m\phi (D_{\eps^{-1}}(\z^{-1}\circ z'))-\DB^m_{z}\phi (D_{\eps^{-1}}(\z^{-1}\circ z))\right|_{z'=e^{\d Yz}}\frac{d\z}{\eps^{Q+2}}\\
&\qquad \quad + \int_{\R^{d+1}}\left|u(e^{\d Y}z)-T_nu(\z,e^{\d Y}z)-\left(u(z)-T_nu(\z,z)\right)\right|
\left|\DB^m_{z}\phi (D_{\eps^{-1}}(\z^{-1}\circ z))\right|\frac{d\z}{\eps^{Q+2}}\\
&\qquad =: \Delta I^{(n,m)}_{\eps,1}+\Delta I^{(n,m)}_{\eps,2}.
\end{align*}
By \eqref{eP1}, Theorem \ref{Taylor1}, Lemma \ref{r2}, and the usual change of variables $\bar \z=D_{\eps^{-1}}(\z^{-1}\circ z)$, we have
\begin{align*}
 \Delta I^{(n,m)}_{\eps,1}&\lesssim \|u\|_{C^{n,\a}_B}\int_{\R^{d+1}}\|\z^{-1}\circ z\|^{n+\a}_B
 \eps^{-m}\left|\DB^m\phi(\bar{\z}')-\DB^m\phi(\bar{\z})\right|_{\bar{\z}'=e^{\frac{\d}{\eps^2}Y}\bar{\z}}\frac{d\z}{\eps^{Q+2}}\\
 &\lesssim \|u\|_{C^{n,\a}_B}\int_{\|\bar\z\|_B\le 1} \eps^{n+\a-m}\left|\frac{\d}{\eps^2}\right|^{\frac 12}\|\phi\|_{C^{m+1,0}_B}d\bar \z
 \lesssim |\d|^{\frac 12}\eps^{n+\a-(m+1)}\|u\|_{C^{n,\a}_B}.
\end{align*}
To control $\Delta I^{(n,m)}_{\eps,2}$ we write $$\Delta I^{(n,m)}_{\eps,2}=\Delta I^{(n,m)}_{\eps,2}\mathbb{I}_{\d>\eps^2}+\Delta I^{(n,m)}_{\eps,2}\mathbb{I}_{\d\le\eps^2}.$$ 
We have
\begin{align*}
\Delta I^{(n,m)}_{\eps,2}\mathbb{I}_{\d>\eps^2}&\le\int_{\R^{d+1}}\left(|u(e^{\d Y}z)-T_nu(\z,e^{\d Y}z)|+|u(z)-T_nu(\z,z)|\right)
\eps^{-m}\left|\DB^m\phi(\bar\z)\right|\frac{d\z}{\eps^{Q+2}}\mathbb{I}_{\d>\eps^2}
\intertext{(using that $\z^{-1}\circ e^{\d Y}z=D_{\eps}e^{\d \eps^{2}Y}\bar \z$ by Lemma \ref{r2} and Remark \ref{r3})}
&\lesssim \|u\|_{C^{n,\a}_B}\int_{\|\bar\z\|_B\le 1} \eps^{n+\a-m}(\|e^{\d\eps^2 Y}\bar \z\|^{n+\a}_B+\|\bar \z\|^{n+\a}_B)
\left|\DB^m\phi(\bar\z)\right|d\bar \z\mathbb{I}_{\d>\eps^2}\\
&\lesssim \eps^{n+\a-m} \|u\|_{C^{n,\a}_B}\mathbb{I}_{\d>\eps^2}\lesssim |\d|^{\frac 12}\eps^{n+\a-(m+1)}\|u\|_{C^{n,\a}_B}.
\end{align*}
Finally, assume for a moment $n\ge 2$. 
By the standard mean value theorem and Lemma \ref{LEM}, there exists $\l\in [0,1]$ such that
\begin{align*}
\Delta I^{(n,m)}_{\eps,2}\mathbb{I}_{\d\le\eps^2}
&\le |\d| \int_{\R^{d+1}}\left|Yu(z')-T_{n-2}Yu(\z,z')\right|_{\z'=e^{\l \d Y}z}
\eps^{-m}\left|\DB^m\phi(\bar\z)\right|\frac{d\z}{\eps^{Q+2}}\mathbb{I}_{\d\le\eps^2}\\
&\lesssim |\d| \|Yu\|_{C^{n-2,\a}_B}\int_{\R^{d+1}}\|\z^{-1}\circ e^{\l\d Y}z\|^{n-2+\a}_B\eps^{-m}
\|\phi\|_{C^{m,0}_B}\frac{d\z}{\eps^{Q+2}}\mathbb{I}_{\d\le\eps^2}\\
&\lesssim |\d| \|u\|_{C^{n,\a}_B}\int_{\|\z\|_B\le 1}\eps^{n-2+\a-m}\|e^{\l\d \eps^2Y}\bar\z\|^{n-2+\a}_B{d\bar \z}\mathbb{I}_{\d\le\eps^2}\\
&\lesssim  |\d|\eps^{n-2+\a-m} \|u\|_{C^{n,\a}_B}\mathbb{I}_{\d\le\eps^2}\lesssim  |\d|^{\frac 12}\eps^{n+\a-(m+1)} \|u\|_{C^{n,\a}_B}.
\end{align*}
If $n=0$ or $n=1$ we have $T_1u(\z,e^{\d Y}z)=T_1u(\z,z)$ and $T_0u(\z,z)=T_0u(\z,e^{\d Y}z)=u(\z)$, so that the expression simplifies and it suffices to exploit the regularity of $u$ w.r.t the field $Y$: this proves \eqref{eP2bis}. 
\medskip\\
\textit{2: Conclusions}
\medskip\\
Let 
$$J^{(n,m)}_{\eps}u(z):=\int_{\R^{d+1}}T_nu(\z,z)\DB^m_z\phi (D_{\eps^{-1}}(\z^{-1}\circ z))\frac{d\z}{\eps^{Q+2}}.$$
By iterative applications of the generalized Leibniz formula and Lemma \ref{LEM}, for $l\le n$:
\begin{align}
\|u-u^{(n)}_{\eps}\|_{C^{l,0}_\a}
&\lesssim \sum_{i=0}^l\left|I^{(n-i,l-i)}_{\eps}\DB^iu\right|_{\infty}+
\sum_{i=0}^{l-1}\left[I^{(n-i,l-i-1)}_{\eps}\DB^iu\right]_{C^{\frac{1}{2}}_Y}\nonumber
\intertext{(using that $\DB^iu\in C^{n-i,\a}_B$ by Theorem \ref{Taylor}, and \eqref{eP2}-\eqref{eP2bis} with $n-i$ instead of $n$ and $m=l-i$)}
&\lesssim\sum_{i=0}^l  \eps^{n+\a-l}\|\DB^iu\|_{C^{n-i,\a}_B}+\sum_{i=0}^{l-1}\ \eps^{n+\a-l}\|\DB^iu\|_{C^{n-i,\a}_B}
\lesssim \eps^{n+\a-l}\|u\|_{C^{n,\a}_B}. \label{eP3}
\end{align}
Similarly, for $l>n$, we have 
\begin{align*}
\|u^{(n)}_{\eps}\|_{C^{l,0}_\a}
\lesssim \sum_{i=0}^n\left(\left|J^{(n-i,l-i)}_{\eps}\DB^iu\right|_{\infty}+
\left[J^{(n-i,l-i-1)}_{\eps}\DB^iu\right]_{C^{\frac{1}{2}}_Y}\right)
\end{align*}
Since $l>n$ any term of the first sum sports, at least, one derivative applied on $\phi$; since the integral of the mollifier \eqref{mollifier} is constant in $z$, the cancellation property $J^{(n-i,l-i)}_{\eps}\DB^iu=I^{(n-i,l-i)}_{\eps}\DB^iu$ holds for any $i=0,\dots, n$. Similarly for the second sum. Therefore we finally get
\begin{align}
\|u^{(n)}_{\eps}\|_{C^{l,0}_\a}
& \lesssim \sum_{i=0}^n\left(\left|I^{(n-i,l-i)}_{\eps}\DB^iu\right|_{\infty}+
\left[I^{(n-i,l-i-1)}_{\eps}\DB^iu\right]_{C^{\frac{1}{2}}_Y}\right)\nonumber\\
& \lesssim\sum_{i=0}^n  \eps^{n+\a-l}\|\DB^iu\|_{C^{n-i,\a}_B}
\lesssim \eps^{n+\a-l}\|u\|_{C^{n,\a}_B}. \label{eP3bis}
\end{align}
The proof is completed.
\endproof

\section{Proofs of Proposition \ref{PROP2} and Theorem \ref{TH}}\label{proof}

\proof[Proof of Proposition \ref{PROP2}]Let $u\in C^{n,0}_B$. By Proposition \ref{PROP} we have 
$$K(\l,u;C^{n_1,0}_B,C^{n_2,0}_B)\lesssim \|u-u_{\eps}\|_{C^{n_1,0}_B}+\l\|u_{\eps}\|_{C^{n_2,0}_B}\lesssim (\eps^{n-n_1}+\l\eps^{n-n_2})\|u\|_{C^{n,0}_B}.$$
Then, taking $\eps=\l^{\frac{1}{n_2-n_1}}$ we directly get $K(\l,u;C^{n_1,0}_B,C^{n_2,0}_B)\lesssim \l^{\frac{n-n_1}{n_2-n_1}}\|u\|_{C^{n,0}_B}$, which gives the set inclusion on the right. It remains to prove that, for any $u\in C^{n_2,0}_B$ 
\begin{equation}\label{eP4}
\|u\|_{C^{n,0}_B}\lesssim \|u\|_{C^{n_1,0}_B}^{1-\frac{n-n_1}{n_2-n_1}}\|u\|_{C^{n_2,0}_B}^{\frac{n-n_1}{n_2-n_1}}.
\end{equation}
We check that \eqref{eP4} holds for $(n_1,n_2)=(n-1,n+1)$, $n\in\N$, in which case, the inequality reads 
\begin{equation}\label{eP4bis}
\|u\|_{C^{n,0}_B}\lesssim \|u\|_{C^{n-1,0}_B}^{\frac{1}{2}}\|u\|_{C^{n+1,0}_B}^{\frac{1}{2}}.
\end{equation}
Then it is straightforward to see that \eqref{eP4} follows by repeated applications of \eqref{eP4bis} with itself. We proceed by induction on $n$. 
By the standard Mean Value Theorem we have
\begin{align*}
|u(e^{\d \p_{x_i}}z)-u(z)-\p_{x_i}u(z)\d|&\le\frac{1}{2}|\p^2_{x_i}u|_{\infty}\d^2, \quad i=1,\dots, p_0\\
|u(e^{\d^2 Y}z)-u(z)|&\le |Yu|_{\infty}\d^2, 
\end{align*}
and therefore 
$$\sum_{i=1}^{p_0}|\p_{x_i}u|_{\infty}+[u]_{C^1_Y}\lesssim \d^{-1}|u|_{\infty}+\d\left(|\p^2_{x_i}u|_{\infty}+|Yu|_{\infty}\right).$$
Taking the optimal $\d>0$ we get 
\begin{equation}\label{eP5}
\|u\|_{C^{1,0}_B}\lesssim |u|_{\infty}+|u|^{\frac 12}_{\infty}\left(|\p^2_{x_i}u|_{\infty}+|Yu|_{\infty}\right)^{\frac 12}\lesssim 
|u|^{\frac 12}_{\infty}\|u\|^{\frac 12}_{C^{2,0}_B}.
\end{equation}
This proves \eqref{eP4bis} for $n=1$. 
Next, by the mean value Theorem along the vector field $Y$ we have that, for every $z\in R^{d+1}$ and $\d\neq 0$, there exists $\bar \d$, $|\bar \d|\le |\d|$ such that 
$$u(e^{\d Y}z)-u(z)-\d Yu(z)=\d \left(Yu(e^{\bar \d Y}z)-Yu(z)\right).$$
Then, dividing by $|\d|^{\frac 12}$ we easily derive
\begin{equation}\label{eP6}
|\d|^{\frac 12}|Yu|_{\infty}\lesssim [u]_{C^{\frac 12}_Y}+|\d|[Yu]_{C^{\frac 12}_Y}
\end{equation}
By \eqref{eP5} and \eqref{eP6}, taking the optimal $\d$, we eventually get 
\begin{align*}
\|u\|_{C^{2,0}_B}&=|u|_{\infty}+|Yu|_{\infty}+\sum_{i=0}^{p_0}\|\p_{x_i}u\|_{C^{1,0}_B}\\
&\lesssim |u|_{\infty}+[u]^{\frac 12}_{C^{\frac 12}_Y}[Yu]^{\frac 12}_{C^{\frac 12}_Y}+\sum_{i=0}^{p_0}|\p_{x_i}u|^{\frac 12}_{\infty}\|\p_{x_i}u\|^{\frac 12}_{C^{2,0}_B}\lesssim \|u\|^{\frac 12}_{C^{1,0}}\|u\|^{\frac 12}_{C^{3,0}}.
\end{align*}
This proves \eqref{eP4bis} for $n=2$. The general case simply follows by the iterative definition of the spaces and the induction hypothesis.
\endproof
Before we proceed with the proof of Theorem \ref{TH}, we recall another important tool in the theory of interpolation, that is the well known Reiteration Theorem (see \cite{MR3753604}, Theorem 1.23 or \cite{Triebel}, Section 1.10.2).
\begin{theorem}[Reiteration Theorem]\label{Reiteration}
Let $\{Z_1,Z_2\}$ be an interpolation couple, and let $E_1,E_2$ be some intermediate spaces between $Z_1$ and $Z_2$. If 
$$\big(Z_1,Z_2\big)_{\theta_i,1}\subseteq E_i\subseteq \big(Z_1,Z_2\big)_{\theta_i,\infty}, \quad i=1,2$$
for some $\theta_i$ such that $0\le \theta_1<\theta_2\le \infty$, then 
$$\big(E_1,E_2\big)_{\a,p}=\big(Z_1,Z_2\big)_{(1-\a)\theta_1+\a\theta_2,p}, \quad \a\in (0,1), \ p\ge 1.$$
\end{theorem}
We are finally ready to prove Theorem \ref{TH}.
\proof[Proof of Theorem \ref{TH}]
The proof is a fairly standard application of Theorem \ref{PROP}, Proposition \ref{PROP2} and Theorem \ref{Reiteration}. We proceed in two steps:
\\
\medskip
\\
\textit{Step 1}: We prove \eqref{eTH} for $\a_1=\a_2=0$ and $n_1=n$, $n_2=n+1$. Precisely, for any $n\in \N_0$, and $0<\alpha<1$, 
$$\left(C^{n,0}_B,C^{n+1,0}_B\right)_{\a,\infty}=\left(C^{n,0}_B,C^{n,1}_B\right)_{\a,\infty}=C^{n,\a}_B.$$
$(\supseteq)$: This inclusion is a direct consequence of the approximation result of Theorem \ref{PROP}. 
Indeed, for any $u\in C^{n,\a}_B$ and $\l\in [0,1)$ we have 
$$K(\l,u;C^{n,0}_B,C^{n+1,0}_B)\lesssim \|u-u_{\eps}\|_{C^{n,0}_B}+\l\|u_{\eps}\|_{C^{n+1,0}_B}\lesssim (\eps^{\a}+\l \eps^{\a-1})\|u\|_{C^{n,\a}_B}.$$
Therefore, taking $\epsilon=\lambda$ we get $$K(\l,u;C^{n,0}_B,C^{n+1,0}_B)\lesssim \l^{\a}\|u\|_{C^{n,\a}_B}.$$  
On the other hand, for $\l\ge 1$ it suffices to take $u_{\eps}\equiv 0$.
\medskip \\
$(\subseteq)$: We proceed by induction on $n$. 
Let $u\in C^{0,\a}_B$, then, for any choice of $a\in C$, $b\in C^{1,0}_B$ such that $u=a+b$ we have $|u|_{\infty}\le |a|_{\infty}+|b|_{\infty}$ and therefore 
$|u|_{\infty}\le K(1,u;C,C^{1,0}_B)\le \|u\|_{(C,C^{1,0}_B)_{\a,\infty}}$.
Similarly, for any $i =1,\dots, p_0$
\begin{equation*}
|u(e^{\delta \p_{x_i}}z)-u(z)|\le 2 |a|_{\infty}+[b]_{C^1_{\p_{x_i}}}|\delta|, \quad 
|u(e^{\delta Y}z)-u(z)|\le 2 |a|_{\infty}+[b]_{C^{\frac 12}_{Y}}|\delta|^{\frac 12}, \quad z\in\R^{d+1},
\end{equation*}
and thus, for any $i =1,\dots, p_0$
\begin{align}
|u(e^{\delta \p_{x_i}}z)-u(z)|&\le 2 K(|\delta|,u;C,C^{1,0}_B)\le 2|\delta|^{\a}\|u\|_{(C,C^{1,0}_B)_{\a,\infty}}, \label{e51}\\
|u(e^{\delta Y}z)-u(z)|&\le 2 K(|\delta|^{\frac 12},u;C,C^{1,0}_B)\le 2|\delta|^{\frac{\a}{2}}\|u\|_{(C,C^{1,0}_B)_{\a,\infty}}\label{e52}.
\end{align}
Gathering together \eqref{e51}-\eqref{e52}, we get 
the inclusion in the case $n=0$.
Next, let $u\in C^{1,\a}_B$. Clearly $|u|_{\infty}\le K(1,u;C^{1,0}_B,C^{2,0}_B)$. 
Moreover,for any choice of $a\in C^{1,0}_B$, $b\in C^{2,0}_B$ such that $u=a+b$ we have, for any $i =1,\dots, p_0$
\begin{align*}
|u(e^{\delta \p_{x_i}}z)-u(z)-\d\p_{x_i}u(z)|&\le 2|\p_{x_i}a|_{\infty}|\d|+\d^2[\p_{x_i}b]_{C^1_{\p_{x_i}}}|\delta^2|, \\ 
|u(e^{\delta Y}z)-u(z)|&\le [a]_{C^{\frac 12}_Y}|\d|^{\frac 12}+|Yb|_{\infty}|\delta|,\\
|\p_{x_i}u(z')-\p_{x_i}u(z)|_{z'=e^{\delta Y}z}&\le 2 |\p_{x_i}a|_{\infty}+[\p_{x_i}b]_{C^{\frac 12}_{Y}}|\delta|^{\frac 12},
\end{align*}
and therefore, for any $i =1,\dots, p_0$
\begin{align}
|u(e^{\delta \p_{x_i}}z)-u(z)-\d\p_{x_i}u(z)|&\le 2|\d|K(|\d|,u;C^{1,0}_B,C^{2,0}_B)\le 2|\d|^{1+\a}\|u\|_{(C^{1,0}_B,C^{2,0}_B)_{\a,\infty}}, \label{e53}\\ 
|u(e^{\delta Y}z)-u(z)|&\le |\d|^{\frac 12}K(|\d|^{\frac 12},u;C^{1,0}_B,C^{2,0}_B)\le |\d|^{\frac{1+\a}{2}}\|u\|_{(C^{1,0}_B,C^{2,0}_B)_{\a,\infty}},\label{e54}\\
|\p_{x_i}u(z')-\p_{x_i}u(z)|_{z'=e^{\delta Y}z}&\le 2 K(|\d|^{\frac 12},u;C^{1,0}_B,C^{2,0}_B)\le 2|\d|^{\frac{\a}{2}}\|u\|_{(C^{1,0}_B,C^{2,0}_B)_{\a,\infty}}. \label{e55}
\end{align}
Gathering together \eqref{e53}-\eqref{e55}, we get 
the inclusion in the case $n=1$. 
Finally, for $u\in C^{n,\a}_B$, $n\ge 2$, we have 
\begin{align*}
\|u\|_{C^{n,\a}_B}&=|u|_{\infty}+\|Yu\|_{C^{n-2,\a}_B}+\sum_{i=1}^{p_0}\|\p_{x_i}u\|_{C^{n-1,\a}_B}\\
&\lesssim  |u|_{\infty}+\|Yu\|_{(C^{n-2,0}_B,C^{n-1,0}_B)_{\a,\infty}}+\sum_{i=1}^{p_0}\|\p_{x_i}u\|_{(C^{n-1,0}_B,C^{n,0}_B)_{\a,\infty}}
\lesssim \|u\|_{(C^{n,0}_B,C^{n+1,0}_B)_{\a,\infty}},
\end{align*}
where we used that $K(\l,\p_{x_i}u;C^{n-1,0}_B,C^{n,0}_B)\le K(\l,u;C^{n,0}_B,C^{n+1,0}_B)$ as well as
 $K(\l,Yu;C^{n-2,0}_B,C^{n-1,0}_B)\le K(\l,u;C^{n,0}_B,C^{n+1,0}_B)$, for any $\l\ge 0$, $i=1,\dots, p_0$.
\\
\medskip
\\
\textit{Step 2}: Let now $n_1, \, n_2 \in \N_0$, $n_1<n_2$, and let $\a_1,\a_2\in [0,1]$. 
By the interpolation result of \textit{Step 1}, and the Reiteration Theorem paired with Proposition \ref{PROP2}, we get
%
\begin{align*}
\left(C^{n_1,\a_1}_B,C^{n_2,\a_2}_B\right)_{\theta,\infty}
&=\left(\left(C^{n_1,0}_B,C^{n_1+1,0}_B\right)_{\a_1,\infty},\left(C^{n_2,0}_B,C^{n_2+1,0}_B\right)_{\a_2,\infty}\right)_{\theta,\infty}\\
&=\left(\left(C^{n_1,0}_B,C^{n_2+1,0}_B\right)_{\frac{\a_1}{n_2+1-n_1},\infty},\left(C^{n_1,0}_B,C^{n_2+1,0}_B\right)_{\frac{n_2+\a_2-n_1}{n_2+1-n_1},\infty}\right)_{\theta,\infty}\\
&=\left(C^{n_1,0}_B,C^{n_2+1,0}_B\right)_{(1-\theta)\frac{\a_1}{n_2+1-n_1}+\theta\frac{n_2+\a_2-n_1}{n_2+1-n_1},\infty}\\
&=:\left(C^{n_1,0}_B,C^{n_2+1,0}_B\right)_{\theta',\infty}
\end{align*}
Taking $n$ as in \eqref{eTHbis} we finally obtain, again by the Reiteration Theorem, 
\begin{align*}
\left(C^{n_1,\a_1}_B,C^{n_2,\a_2}_B\right)_{\theta,\infty}&=\left(C^{n,0}_B,C^{n+1,0}_B\right)_{\theta'(n_2+1-n_1)-(n-n_1),\infty}
=C^{n,\a_1+\theta[(n_2+\a_2)-(n_1+\a_1)]-(n-n_1)}_B.
\end{align*}
The proof is completed. 
\endproof

\section*{Acknowledgments}
This research was supported by the Gruppo Nazionale per l'Analisi Matematica, la Probabilit\`a e le loro Applicazioni (GNAMPA) of the Istituto Nazionale di Alta Matematica (INdAM), and the University of Bologna. We would like to thank Professor A. Pascucci for his interest in our work and for several helpful conversations.

\bibliographystyle{acm}
\bibliography{bib1}

\def\cprime{$'$} \def\cprime{$'$} \def\cprime{$'$}
  \def\lfhook#1{\setbox0=\hbox{#1}{\ooalign{\hidewidth
  \lower1.5ex\hbox{'}\hidewidth\crcr\unhbox0}}} \def\cprime{$'$}
  \def\cprime{$'$} \def\cprime{$'$} \def\cprime{$'$} \def\cprime{$'$}
  \def\polhk#1{\setbox0=\hbox{#1}{\ooalign{\hidewidth
  \lower1.5ex\hbox{`}\hidewidth\crcr\unhbox0}}}
\begin{thebibliography}{10}

\bibitem{Adams1975}
{\sc Adams, R.~A.}
\newblock {\em Sobolev spaces}.
\newblock Academic Press, New York-London, 1975.
\newblock Pure and Applied Mathematics, Vol. 65.

\bibitem{BarucciPolidoroVespri}
{\sc Barucci, E., Polidoro, S., and Vespri, V.}
\newblock Some results on partial differential equations and {A}sian options.
\newblock {\em Math. Models Methods Appl. Sci. 11}, 3 (2001), 475--497.

\bibitem{MR0230022}
{\sc Butzer, P.~L., and Berens, H.}
\newblock {\em Semi-groups of operators and approximation}.
\newblock Die Grundlehren der mathematischen Wissenschaften, Band 145.
  Springer-Verlag New York, Inc., New York, 1967.

\bibitem{Cercignani}
{\sc Cercignani, C.}
\newblock {\em The {B}oltzmann equation and its applications}.
\newblock Springer-Verlag, New York, 1988.

\bibitem{MR4358660}
{\sc Chaudru~de Raynal, P.-E., and Menozzi, S.}
\newblock Regularization effects of a noise propagating through a chain of
  differential equations: an almost sharp result.
\newblock {\em Trans. Amer. Math. Soc. 375}, 1 (2022), 1--45.

\bibitem{CdRMPZ22}
{\sc Chaudru~de Raynal, P.-E., Menozzi, S., Pesce, A., and Xicheng, Z.}
\newblock Heat kernel and gradient estimates for kinetic {SDE}s with low
  regularity coefficients.
\newblock {\em arXiv:2203.11515\/} (2022).

\bibitem{CIBELLI201987}
{\sc Cibelli, G., Polidoro, S., and Rossi, F.}
\newblock Sharp estimates for {G}eman-{Y}or processes and applications to
  {A}rithmetic {A}verage {A}sian options.
\newblock {\em Journal de Math\'ematiques Pures et Appliqu\'ees 129\/} (2019),
  87--130.

\bibitem{Desvillettes}
{\sc Desvillettes, L., and Villani, C.}
\newblock On the trend to global equilibrium in spatially inhomogeneous
  entropy-dissipating systems: the linear {F}okker-{P}lanck equation.
\newblock {\em Comm. Pure Appl. Math. 54}, 1 (2001), 1--42.

\bibitem{DiFrancescoPascucci2}
{\sc Di~Francesco, M., and Pascucci, A.}
\newblock On a class of degenerate parabolic equations of {K}olmogorov type.
\newblock {\em AMRX Appl. Math. Res. Express 3\/} (2005), 77--116.

\bibitem{DiFrancescoPolidoro}
{\sc Di~Francesco, M., and Polidoro, S.}
\newblock Schauder estimates, {H}arnack inequality and {G}aussian lower bound
  for {K}olmogorov-type operators in non-divergence form.
\newblock {\em Adv. Differential Equations 11}, 11 (2006), 1261--1320.

\bibitem{FrentzNystromPascucci2010}
{\sc Frentz, M., Nystr{\"o}m, K., Pascucci, A., and Polidoro, S.}
\newblock Optimal regularity in the obstacle problem for {K}olmogorov operators
  related to {A}merican {A}sian options.
\newblock {\em Math. Ann. 347}, 4 (2010), 805--838.

\bibitem{Imbert}
{\sc Imbert, C., and Mouhot, C.}
\newblock The {S}chauder estimate in kinetic theory with application to a toy
  non-linear model.
\newblock {\em Annales H. Lebesgue 4\/} (2021), 369--405.

\bibitem{Imbert2}
{\sc Imbert, C., and Silvestre, L.}
\newblock The {S}chauder estimate for kinetic integral equations.
\newblock {\em Analysis and PDE 14\/} (2021), 171--204.

\bibitem{MR1406091}
{\sc Krylov, N.~V.}
\newblock {\em Lectures on elliptic and parabolic equations in {H}\"{o}lder
  spaces}, vol.~12 of {\em Graduate Studies in Mathematics}.
\newblock American Mathematical Society, Providence, RI, 1996.

\bibitem{LanconelliPolidoro}
{\sc Lanconelli, E., and Polidoro, S.}
\newblock On a class of hypoelliptic evolution operators.
\newblock {\em Rend. Sem. Mat. Univ. Politec. Torino 52}, 1 (1994), 29--63.

\bibitem{Lucertini}
{\sc Lucertini, G., Pagliarani, S., and Pascucci, A.}
\newblock Optimal regularity for degenerate {K}olmogorov equations with rough
  coefficients.
\newblock {\em To appear in J. Evol. Equ.\/} (2023).

\bibitem{Lucertini2}
{\sc Lucertini, G., Pagliarani, S., and Pascucci, A.}
\newblock Optimal {S}chauder estimates for kinetic {K}olmogorov equations with
  time measurable coefficients.
\newblock {\em arXiv:2304.13392\/} (2023).

\bibitem{Lunardi2}
{\sc Lunardi, A.}
\newblock {\em Analytic {S}emigroups and {O}ptimal {R}egularity in {P}arabolic
  {P}roblems}.
\newblock Modern Birkh\"{a}user Classics. Birkh\"{a}user Basel, 1995.

\bibitem{Lunardi}
{\sc Lunardi, A.}
\newblock Schauder estimates for a class of degenerate elliptic and parabolic
  operators with unbounded coefficients in ${\R}^{N}$.
\newblock {\em Ann. Scuola Norm. Sup. Pisa Cl. Sci. (4) 24}, 1 (1997),
  133--164.

\bibitem{MR3753604}
{\sc Lunardi, A.}
\newblock {\em Interpolation theory}, vol.~16 of {\em Lecture Notes. Scuola
  Normale Superiore di Pisa}.
\newblock Edizioni della Normale, Pisa, 2018.

\bibitem{MR1751429}
{\sc Manfredini, M.}
\newblock The {D}irichlet problem for a class of ultraparabolic equations.
\newblock {\em Adv. Differential Equations 2}, 5 (1997), 831--866.

\bibitem{Manfredini}
{\sc Manfredini, M.}
\newblock The {D}irichlet problem for a class of ultraparabolic equations.
\newblock {\em Adv. Differential Equations 2}, 5 (1997), 831--866.

\bibitem{MR3429628}
{\sc Pagliarani, S., Pascucci, A., and Pignotti, M.}
\newblock Intrinsic {T}aylor formula for {K}olmogorov-type homogeneous groups.
\newblock {\em J. Math. Anal. Appl. 435}, 2 (2016), 1054--1087.

\bibitem{MR4355925}
{\sc Pascucci, A., and Pesce, A.}
\newblock On stochastic {L}angevin and {F}okker-{P}lanck equations: the
  two-dimensional case.
\newblock {\em J. Differential Equations 310\/} (2022), 443--483.

\bibitem{PP22}
{\sc Pascucci, A., and Pesce, A.}
\newblock Sobolev embeddings for kinetic {F}okker-{P}lanck equations.
\newblock {\em arXiv:2209.05124\/} (2022).

\bibitem{MR2328004}
{\sc Tartar, L.}
\newblock {\em An introduction to {S}obolev spaces and interpolation spaces},
  vol.~3 of {\em Lecture Notes of the Unione Matematica Italiana}.
\newblock Springer, Berlin; UMI, Bologna, 2007.

\bibitem{Triebel}
{\sc Triebel, H.}
\newblock {\em Interpolation {T}heory, {F}unction {S}paces, {D}ifferential
  {O}perators}, vol.~18 of {\em North-Holland Mathematical Library}.
\newblock 1978.

\end{thebibliography}

\end{document}